# Solution of variable order fractional differential equations using Homotopy Analysis Method

**Vivek Mishra, S. Das**

**Abstract** In the present article an endeavor is made to solve the variable order fractional diffusion equations using a powerful method viz., Homotopy Analysis method. It is demonstrated how the method can be used while solving approximately two types of variable order fractional diffusion equations having physical importance. Numerical simulation results show that the method is reliable and effective for solving fractional order diffusion equations even when the order of the derivative is varying with respect to space or time or both or it is dependent upon some other parameters.

**Keywords :** Deformation Equation, Convergence control parameter, Residual error, series solution

**2010 Mathematics Subject Classification**: 26A33, 34G20, 35A20, 34A08, 65Mxx

## 1. Introduction

Fractional calculus ([1]-[4]) gives us flexibility to do the integration and differentiation of arbitrary order. Till date constant order fractional derivative and integration are tremendously used in the modeling of the physical problem of fluid flow or diffusion .The flow of the fluid in heterogeneous medium is becoming an important area of research in recent time. In homogeneous medium the diffusion equation of constant order explains all the stages of the diffusion like normal diffusion, anomalous diffusion etc. In heterogeneous media like porous media where different permeability exists at different places so to explain the diffusion process in this kind of media, variable order fractional diffusion equation is required. In many diffusion processes, diffusion speed changes with time either it is accelerating or it is retarding. Traditionally researchers have tried to model this phenomenon of changing speed by taking time dependent diffusion coefficient. The variable order fractional differential equation is the appropriate selection for the modeling of this kind of phenomenon. If the diffusion process is interfered by any kind of disturbances or by noises due to which there is fluctuation in the whole system then the system will show deviation from constant order fractional system to variable order system.

In physical process, memory also plays a very important role and the best way to use the effect of the memory on the process is the application of constant order fractional system but if the memory property of the system changes with time or with space or with both then variable order partial differential equations provide a better description of the memory. Recently, many researchers have found that variable order model is the most appropriate for the modeling of the physical process where behavior of fractional order derivative changes with time and space.

Caputo [5-6] first gave the concept of distributed order derivative in 1967. Distributed order derivative was not sufficient for the modeling of diffusion process with varying speed. Samko et al. [7-8] first proposed the concept of variable order operator. Lorenzo and Hartley [9] have investigated the different types of variable order fractional order operator definitions. Smit and Vries [10] investigated stress strain behavior of viscoelastic material using fractional order differential equation and shown that order '$\alpha$' depends on strain level. Sun [11] used variable order model for describing the different kinds of anomalous diffusion process. Coimbra [12] has discussed the mechanics of variable order differential operator. Metzler et al. [13] have explained that order of fractional differential equation used in relaxation process and reaction dynamics depends on temperature. Kobelev et al. [14] have studied statistical and dynamical systems with fixed and variable memories. Cooper et al. [15] have used variable order model for processing the geographical data verification. Lin et al. [16] discussed the stability and convergence of an explicit finite difference approximation for variable order nonlinear diffusion equation. Pedro et al. [17] investigated the motion of the particle suspended in a viscous fluid with drag force using variable order calculus. Sun et al.[18] have done a comparative study of constant order and variable order fractional models. Razminia et al.[19] have discussed the existence for non- autonomous variable order fractional differential equations. Sun et al.[20] have made the study of mean square displacement behaviors of anomalous diffusion with variable and random orders.

Many researchers have successfully solved the variable order fractional differential equation using numerical methods. Valerio et al. [21] have given the numerical approximation of variable order fractional derivatives. Zhang et al.[22] have solved the variable fractional order mobile immobile advection dispersion model. Chen et al. [23] have given numerical simulation of a new two dimensional variable order fractional percolation equation in non homogeneous porous media. Shen [24] has given Numerical technique for the variable order time fractional diffusion equation. Sun et al. [25] have explained a finite difference scheme for variable order time fractional diffusion equation. Sun et al. [26] have used variable indexed model to explain transient dispersion in heterogeneous media. To the best of authors' knowledge no researcher has tried to solve the variable order fractional differential equation using Homotopy analysis method.

Homotopy analysis method (HAM) is proposed by Liao [27] in 1992, where a simple concept of topology is used known as Homotopy. Since then Liao himself has given strong effort to prove the applicability and validity of HAM ([28]-[35]). The method is applied by many researchers ([36]-[38]) in many fields to solve the Mathematical models of integer order as well as fractional order. The basic concept in HAM is to give the successive approximate solutions of the problem which leads to the exact solution of the problem. The benefit of this method is that it is free from small /large physical parameters. This method also ensures the convergence of the series solution and provides a great flexibility to choose the initial approximation for the solution of the problem. The beauty of the present article is successful implementation of HAM to solve the variable order fractional diffusion equations where the convergence of the series form of the solution is obtained by minimizing the averaged residual error.

## 2. Formula and definition

**2.1** The variable order fractional derivative of Caputo type is defined by Coimbra [12] as

$$^C D_t^{\alpha(t)} f(t) = \frac{1}{\Gamma(1-\alpha(t))} \int_{0^+}^{t} \frac{f'(\tau)}{(t-\tau)^{\alpha(t)}} d\tau + \frac{f(0_+) - f(0_-) \, t^{\alpha(t)}}{\Gamma(1-\alpha(t))}, \qquad 0 < \alpha(t) < 1. \tag{1}$$

For simplicity it is assumed that the function is well defined at $t=0$ and hence the definition takes the structure as follows

$$^C D_t^{\alpha(t)} f(t) = \frac{1}{\Gamma(1-\alpha(t))} \int_{0^+}^{t} \frac{f'(\tau)}{(t-\tau)^{\alpha(t)}} d\tau, \qquad 0 < \alpha(t) < 1. \tag{2}$$

**2.2** Samko [7] has proposed the definition of variable order integration as

$$J_t^{\alpha(t)} f(t) = \frac{1}{\Gamma(\alpha(t))} \int_{0^+}^{t} (t-\tau)^{\alpha(t)-1} f(\tau) \, d\tau, \qquad \text{Re}(\alpha(t)) > 0. \tag{3}$$

Using the above definitions, the following formulas are obtained (Sun et al. [20])

$$^C D_t^{\alpha(t)} t^\beta = \frac{\Gamma(1+\beta)}{\Gamma(\beta - \alpha(t) + 1)} t^{\beta - \alpha(t)}, \qquad 0 < \alpha(t) < 1 \tag{4}`$$

and

$$J_0^{\alpha(t)} t^\beta = \frac{\Gamma(1+\beta)}{\Gamma(\beta + \alpha(t) + 1)} t^{\beta + \alpha(t)}, \qquad 0 < \alpha(t) < 1. \tag{5}$$

## 3. Solution of the problems by HAM

**3.1** We consider the following variable order time fractional Diffusion equation without source term which is already solved numerically by Sun et al. [25]:

$$^c_0 D_t^{\alpha(x,t)} u(x,t) = K \frac{\partial^2 u(x,t)}{\partial x^2}, \tag{6}$$

where $u(x,0) = \sin(\frac{\pi x}{L})$, with $u(0,t) = u(L,t) = 0$,

$$\alpha(x,t) = 0.8 + \frac{0.2xt}{LT} \tag{7}$$

Considering the Linear operator as

$$L[u(x,t)] = {}^c_0 D_t^{\alpha(x,t)} u(x,t) \tag{8}$$

and the nonlinear operator as

$$N[u(x,t)] = {}_0^c D_t^{\alpha(x,t)} u(x,t) - K\frac{\partial^2 u}{\partial x^2} \qquad (9)$$

The m-th order deformation equation is

$$u_m(x,t) = \chi_m u_{m-1}(x,t)] + \hbar J_t^{\alpha(x,t)} R_m[\vec{u}_{m-1}(x,t)] + c , \qquad (10)$$

where $\chi_m = \begin{cases} 0 & m \leq 1 \\ 1 & m > 1 \end{cases}$ \qquad (11)

and c can be calculated using the initial conditions.

Taking $u_0(x,t) = \sin(\frac{\pi x}{L})$, we get \qquad (12)

$$u_1(x,t) = \frac{\hbar \pi^2 k}{L^2} \sin(\frac{\pi x}{L}) \frac{t^{\alpha(x,t)}}{\Gamma(1+\alpha(x,t))} \qquad (13)$$

$$u_2(x,t) = \frac{\hbar (\hbar+1)\pi^2 k}{L^2} \sin(\frac{\pi x}{L}) \frac{t^{\alpha(x,t)}}{\Gamma(1+\alpha(x,t))} + \frac{\hbar^2 \pi^4 k^2}{L^4} \sin(\frac{\pi x}{L}) \frac{t^{2\alpha(x,t)}}{\Gamma(1+2\alpha(x,t))} \qquad (14)$$

$$u_3(x,t) = \frac{\hbar (\hbar+1)^2 \pi^2 k}{L^2} \sin(\frac{\pi x}{L}) \frac{t^{\alpha(x,t)}}{\Gamma(1+\alpha(x,t))} + \frac{2\hbar^2 (\hbar+1)\pi^4 k^2}{L^4} \sin(\frac{\pi x}{L}) \frac{t^{2\alpha(x,t)}}{\Gamma(1+2\alpha(x,t))}$$

$$+ \frac{\hbar^3 \pi^6 k^3}{L^6} \sin(\frac{\pi x}{L}) \frac{t^{3\alpha(x,t)}}{\Gamma(1+3\alpha(x,t))} \qquad (15)$$

$$u_4(x,t) = \frac{\hbar (\hbar+1)^3 \pi^2 k}{L^2} \sin(\frac{\pi x}{L}) \frac{t^{\alpha(x,t)}}{\Gamma(1+\alpha(x,t))} + \frac{3\hbar^2 (\hbar+1)^2 \pi^4 k^2}{L^4} \sin(\frac{\pi x}{L}) \frac{t^{2\alpha(x,t)}}{\Gamma(1+2\alpha(x,t))}$$

$$+ \frac{3\hbar^3 (\hbar+1)\pi^6 k^3}{L^6} \sin(\frac{\pi x}{L}) \frac{t^{3\alpha(x,t)}}{\Gamma(1+3\alpha(x,t))} + \frac{\hbar^3 \pi^8 k^4}{L^8} \sin(\frac{\pi x}{L}) \frac{t^{4\alpha(x,t)}}{\Gamma(1+4\alpha(x,t))} \qquad (16)$$

Proceeding in the similar manner, we can calculate the other terms $u_n$, $n \geq 5$ and hence finally the approximate series solution of the problem is obtained as

$$u(x,t) = \lim_{N \to \infty} \phi_N(x,t), \qquad (17)$$

where $\phi_N(t) = \sum_{n=0}^{N-1} u_n(t)$, $N \geq 1$.

**3.2** Consider the Diffusion equation as

$$D_0^{\alpha(x,t)} u(x,t) = \frac{\partial}{\partial x}\left(u(x,t)\frac{\partial}{\partial x} u(x,t)\right) + u(1-u) , \tag{18}$$

where $\alpha(x,t) = xt$ and $u(x,0) = x$, u(0, t) = u(1, t) = 0.

Considering the linear operator as

$$L[u(x,t)] = {}_0^c D_t^{\alpha(x,t)} u(x,t) \tag{19}$$

and the nonlinear operator as

$$N[u(x,t)] = {}_0^c D_t^{\alpha(x,t)} u(x,t) - \frac{\partial}{\partial x} u(x,t) \frac{\partial}{\partial x} u(x,t) - u(x,t) \frac{\partial^2}{\partial x^2} u(x,t) - u(x,t) + u^2(x,t) , \tag{20}$$

the m-th order deformation equation is obtained as

$$u_m(x,t) = \chi_m u_{m-1}(x,t) + \hbar J_t^{\alpha(x,t)} R_m[\vec{u}_{m-1}(x,t)] + c , \tag{21}$$

where $\chi_m$ is defined in (11).

Taking $u_0 = x$, we obtain

$$u_1 = \hbar(x^2 - x - 1)\frac{t^{\alpha(x,t)}}{\Gamma(1+\alpha(x,t))} \tag{22}$$

$$u_2 = \hbar(\hbar+1)(x^2 - x - 1)\frac{t^{\alpha(x,t)}}{\Gamma(1+\alpha(x,t))} + \hbar^2(2x^3 - 3x^2 - 7x + 3)\frac{t^{2\alpha(x,t)}}{\Gamma(1+2\alpha(x,t))} \tag{23}$$

$$u_3 = \hbar(\hbar+1)^2 (x^2 - x - 1)\frac{t^{\alpha(x,t)}}{\Gamma(1+\alpha(x,t))} + \hbar^2(\hbar+1)(4x^3 - 6x^2 - 14x + 6)\frac{t^{2\alpha(x,t)}}{\Gamma(1+2\alpha(x,t))}$$

$$+ \hbar^3(4x^4 - 8x^3 - 35x^2 + 31x + 11)\frac{t^{3\alpha(x,t)}}{\Gamma(1+3\alpha(x,t))}$$

$$+ \hbar^3(x^4 - 2x^3 - 7x^2 + 8x + 2)\frac{\Gamma(1+2\alpha(x,t))t^{3\alpha(x,t)}}{(\Gamma(1+\alpha(x,t)))^2 \Gamma(1+3\alpha(x,t))} \tag{24}$$

Finally the approximate series solution of the problem is obtained using equation (17).

In this article the convergence control parameter of the HAM is evaluated by reducing the residual error as given by Liao [38] at the m-th order approximation. The exact square residual error is given by

$$E_m = \iint_\Omega \left( N\left[ \sum_{i=0}^m u_i(x,t) \right] \right)^2 dx\, dt \tag{25}$$

However to calculate the residual error and to reduce the CPU time, the authors have used the method suggested by Liao[38] to calculate residual error as

$$E_m = \frac{1}{(M_x+1)(M_t+1)} \sum_{J=1}^{M_x} \sum_{k=1}^{M_t} \left( N\left[ \sum_{i=0}^{m} u_i(j\Delta x., k\Delta t) \right] \right)^2 , \qquad (26)$$

where the optimal value of $\hbar$ is calculated by minimizing the averaged residual error through the formula

$$\frac{\partial E_m}{\partial \hbar} = 0 . \qquad (27)$$

## 4. Results and Discussion

In this section during the analysis of the residual error of the first considered problem, it is assumed that $M_x=34$, $M_t=34$ and $\Delta x = \frac{10}{34}, \Delta t = \frac{1}{34}$ and the obtained results are displayed through Table 1. It is seen from the table that the $E_m$ decreases if the number of term in the series solution increases. Variation of $E_m$ for various $\hbar$ is shown through Fig.1 for k=0.01, $\alpha = 0.8 + 0.2\frac{xt}{LT}, L=1, T=10$.

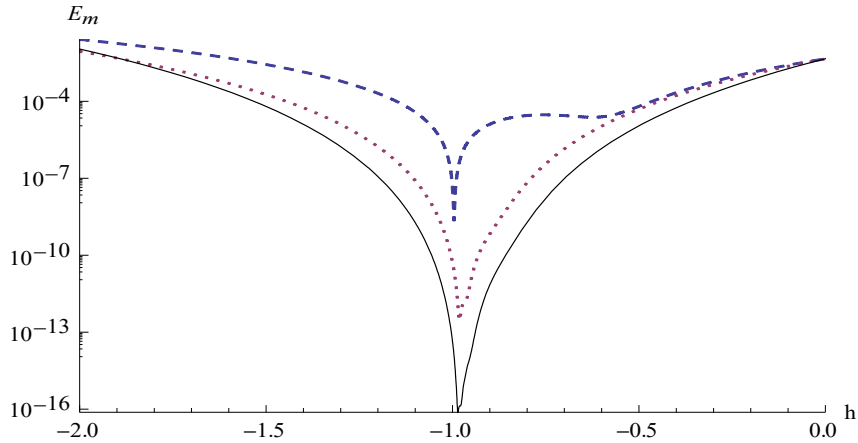

**Fig.1.** Plots of exact residual error $E_m$ versus $\hbar$.

| No of terms in the series solution | Residual error ($E_m$) | Value of $\hbar$ |
|---|---|---|
| Three terms | 2.13959×10$^{-9}$ | -0.995985 |
| Four terms | 9.33408×10$^{-13}$ | -0.975296 |
| Five terms | 2.1842×10$^{-16}$ | -0.975296 |

**Table.1.** comparison of residual error ($E_m$) and convergence control parameter ($\hbar$)

It is also seen from the Fig.1 and Table.1 that the magnitude of residual error is minimum at $\hbar = -0.975296$, when five terms in the series solution are taken. For this value of $\hbar$ the variation of $u(x,t)$ vs. time is depicted through Fig.2 which is in complete agreement with the results given in [25]

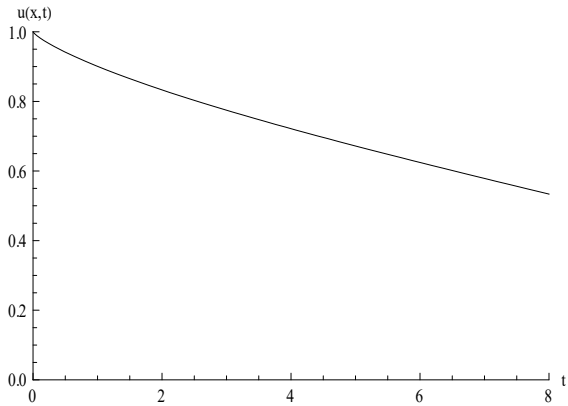

**Fig.2.** Evolution curve of $u(x=0.5,t)$ vs. t at $\hbar = -0.975296$.

Similar calculations are done for the second considered problem to achieve the value of $\hbar$ for which the residual error is minimum. Taking $M_x = 10$, $M_t = 10$ and $\Delta x = \frac{1}{10}, \Delta t = \frac{1}{10}$ the corresponding results are displayed through Table.2 and Fig.3. The variation of $u(x,t)$ for $\hbar = -0.13456$ are displayed through Fig.4 for different values of x and t as and when $\alpha = xt$.

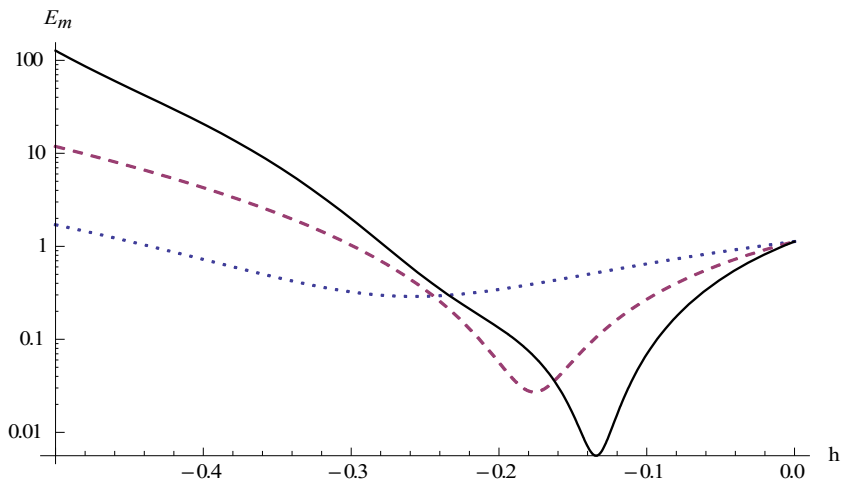

**Fig.3.** Plots of exact residual error $E_m$ versus $\hbar$ for $\alpha = xt$.

| No of terms in the series solution | Residual error ($E_m$) | Value of $\hbar$ |
|---|---|---|
| Two terms | 0.289179 | -0.257313 |
| Three terms | 0.027088 | -0.176075 |
| Four terms | 0.00560894 | -0.134256 |

**Table.2.** Comparison of residual error ($E_m$) and convergence control parameter ($\hbar$)

taking the $\hbar = -0.134256$ we have plotted the solution curve of the above problem

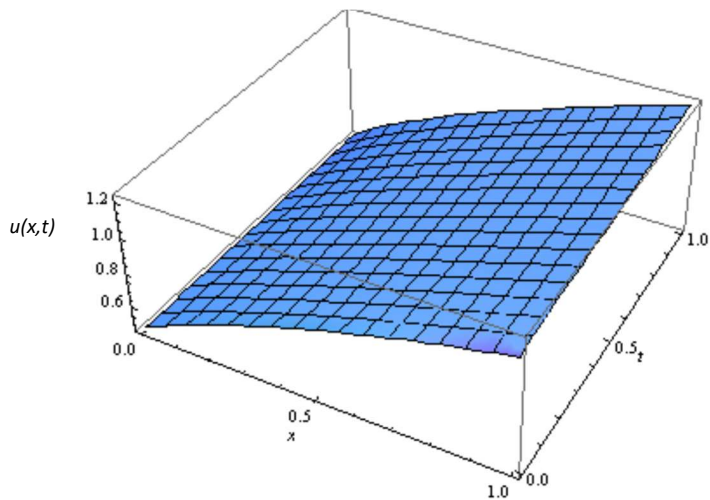

**Fig.4.** Evolution curve of u(x, t) vs x and t.

Solution of the first and second problem is displayed through Fig.2. and Fig.4.

## 5. Conclusion

In the present article authors' first aim was to solve the variable order Fractional differential equation and matching the result with [25] taking the same values of parameters to authenticate the effectiveness of the method even for variable order problems. Our second aim was to solve a nonlinear variable order fractional diffusion equation with a reaction term and furnish the stability of the solution for the different particular cases which is first of its kind in case of a variable order fractional diffusion system in nonlinear dynamics.